\documentclass{article}
\usepackage{color}
\makeatletter \@addtoreset{equation}{subsection}

\makeatother \makeatletter

\usepackage{amsmath,amsfonts,calrsfs,amssymb,color,verbatim,eucal,mathrsfs}

\newtheorem{theorem}{\bf Theorem}[section]

\newcommand{\C}{{\mathbb{C}}}
\newcommand{\N}{{\mathbb{N}}}

\newcommand{\R}{\mathbb{R}}

\def\borel{\mathcal B}

\title{Ornstein-Uhlenbeck semigroups in infinite dimension}

\author{A. Lunardi\thanks{Dipartimento di Scienze Matematiche, Fisiche e Informatiche, 
Universit\`a di Parma, Parco Area delle Scienze 53/A, 43124 Parma, Italy; 
e-mail alessandra.lunardi@unipr.it}  \, and 
D. Pallara\thanks{Dipartimento di Matematica e Fisica "Ennio De Giorgi", Universit\`a del Salento, and INFN, 
Sezione di Lecce. POB 193, 73100 Lecce, Italy; e-mail diego.pallara@unisalento.it}}
\date{}
\begin{document}

\maketitle

\begin{abstract}
This is a survey paper about Ornstein-Uhlenbeck semigroups in infinite dimension, and their generators. 
We start from the classical Ornstein-Uhlenbeck semigroup in Wiener spaces and then discuss the general 
case in Hilbert spaces. Finally, we present some results for O-U semigroups in Banach spaces. 
\end{abstract}


\section*{Introduction}
In this article we present the basic results on Ornstein-Uhlenbeck semigroups in infinite dimensional spaces.  

After an introductory section with notation and preliminaries, the classical O-U semigroup in separable Banach spaces is  discussed in Section \ref{sect:OUclassico}; we refer to the  survey paper \cite{BogaSurvey} for many  details and historical notes. 

The main body of this paper is Section  \ref{sect:Hilbert}, where we describe the theory of Ornstein-Uhlenbeck semigroups in separable Hilbert spaces. We  refer to the book \cite{DPZbrutto} for detailed descriptions of the basic ideas in simple cases, examples and applications, although more precise and more general results (as well as subsequent developments) are elsewhere. 

In the last section we consider Ornstein-Uhlenbeck semigroups in separable Banach spaces. 
There are many more technicalities and many less examples than in the Hilbert setting, and in this short survey we have not room to give details, so we only  briefly list some extensions of the results of Section  \ref{sect:Hilbert} to the  Banach case.

\section{Preliminaries}\label{SecPrelimin}

Throughout the paper  $X$ is a separable real Banach space, with norm $\|\cdot\|$. $\mathcal B_b(X)$ (resp. $C_b(X)$, 
$BUC(X)$) is the space of the Borel measurable (resp. continuous, uniformly continuous) and bounded functions from $X$ to 
$\R$, endowed with the sup norm $\|\cdot\|_{\infty}$. Occasionally, we will be concerned also with the 
{\em mixed topology} in $C_b(X)$, 
which is the finest locally convex topology that agrees with the topology of uniform convergence on compact sets on every bounded set in $X$. 

If $Y$ is any Banach space, $\mathscr L(X,Y)$ is the space of the linear bounded operators from $X$ to $Y$; as usual if $Y=X$ it is denoted by $\mathscr L(X)$ and if $X=\R$ it is denoted by 
$X^*$. For $2\leq h\in \N$,  $\mathcal L^h(X)$ is the space of the continuous $h$-linear operators from $X^h$ to $\R$. 

The Borel $\sigma$-algebra  $\borel(X)$ coincides with the $\sigma$-algebra ${\mathscr E}(X)$ generated by 
the {\em cylindrical} sets, i.e, the sets of the form
$C=\{x\in X:\ \left(f_1(x),\ldots,f_n(x)\right)\in C_0\}$,
where $f_1,\ldots,f_n\in X^*$ and $C_0\in \borel(\R^n)$, see e.g.  \cite[Ch. 1]{VTC}. Accordingly, a function $f:X\to\R$ is called {\em cylindrical} if there are $f_1,\ldots,f_n\in X^*$ and $\varphi:\R^n\to\R$ such that $f(x)=\varphi(f_1(x),\ldots,f_n(x))$. 

If $X$ is a Hilbert space, we denote by $\langle \cdot, \cdot\rangle$ its inner product, and we denote 
by $\mathscr L_1(X)$, $\mathscr L_2(X)$   the subspaces of ${\mathscr L}(X)$ consisting  of the nuclear self-adjoint operators, and by the Hilbert-Schmidt operators, respectively.  

\subsection{Symmetric and positive operators.} 
An operator  $Q\in \mathscr L(X^*, X)$ is called {\em symmetric} if $g(Qf) = f(Qg)$ for every $f$, $g\in X^*$, and {\em positive} if $f(Qf)\geq 0$ for every $f\in X^*$ (in fact, the right word should be ``nonnegative" but we adopt the common terminology). As usual, if $X$ is a Hilbert space we identify $X$ and $X^*$, and the above notions correspond respectively to  a self-adjoint and nonnegative $Q\in  \mathscr L(X)$. 

If $Q$ is symmetric and positive, there exists a unique Hilbert space $H_Q$ continuously embedded in $X$ 
such that $Q(X^*)$ is dense in $H_Q$ and $\langle Qf, Qg\rangle _{H_Q} = g(Qf)$, for every $f$, $g\in X^*$ 
(\cite[Prop. III.1.6]{VTC}). Denoting by $i:  H_Q\to X$ the embedding we have $\|i\|_{\mathscr L(H_Q, X)} = \|Q\|_{\mathscr L(X^*, X)}^{1/2}$ and 
$i\circ i^* = Q$; see \cite[Chapter III]{VTC}. 
$H_Q$ may be equivalently constructed by completing  $Q(X^*)$ with respect to  the  norm associated with the inner product $(Qf, Qg)\mapsto g(Qf)$. It is easily seen that every Cauchy sequence $(Qf_n)$ 
in such norm converges in $X$, and two equivalent Cauchy sequences converge in $X$ to the same limit. Identifying (the 
equivalence class of) any Cauchy sequence $(Qf_n)$ with its $X$-limit $h$,  the completion is identified with  $H_Q$. 

If $X$ is a Hilbert space and $Q\in  \mathscr L(X)$ is self-adjoint and nonnegative, $H_Q$ is just $Q^{1/2}(X)$ with the  inner product  $\langle Q^{1/2}x,  Q^{1/2}y\rangle_{H_Q}  = \langle x, y\rangle$ for every $x$, $y\in X$, or equivalently 
$\langle h, k\rangle_{H_Q} = \langle Q^{-1/2}h, Q^{-1/2}k\rangle$. Here, if $Q^{1/2}$  is not one to one,   $Q^{-1/2}$ denotes its pseudo-inverse$^($\footnote{If   $T\in \mathscr L(X)$ is self-adjoint and nonnegative, for every $h\in T (X)$ 
$T^{-1}h$ is the element of minimal norm in the set $T^{-1 }(\{h\})$. We have $T^{-1 }h = Py$ for every $y\in T(\{h\})$, where  $P$ is the orthogonal projection on $\overline{T (X)} = ($ Ker$\,T)^{\perp}$.}$^)$.

The space $H_Q$ is sometimes called {\it reproducing kernel Hilbert space} associated with $Q$, but since the expression ``reproducing kernel Hilbert space" has several different meanings in the literature, we will not use it. 

\subsection{Regular functions.}

Let $Y$ be any Banach space,    $\alpha \in (0,1)$, $k\in \N$. 

$C_b^{\alpha}(X;Y)$ is the space of the bounded and $\alpha$-H\"older continuous functions from $X$ to $Y$, endowed with its standard norm $\|f\|_{C^\alpha_b(X;Y)} := \|f\|_{\infty} + [f]_{C^\alpha(X;Y)}$, where 
$ [f]_{C^\alpha(X;Y)} = \sup_{x, \, y\in X; \; x\neq y} \|f(x) - f(y)\|_Y/\|x-y\|^{\alpha}$. 
If $Y=\R$ we set $C_b^{\alpha}(X;\R)=: C_b^{\alpha}(X)$. 
 
$C^k_b(X)$ is the space  of the $k$ times Fr\'echet differentiable functions $F:X\to \R$, with continuous and bounded derivatives $D^jf : X\to \mathcal L^j(X)$ for every $j=1, \ldots, k$. 
The first order Fr\'echet derivative $D^1$ is  denoted by $D$. 

If $X$ is a Hilbert space and $f:X\to \R$ is Fr\'echet differentiable at $x$,  by the Riesz isometry there is a unique 
$z\in X$ such that $Df(x)(h) =\langle z, h\rangle $ for every $h\in X$. Such $z$ is denoted by $\nabla f(x)$. 

$C_b^{k+\alpha}(X)$ is the space  of the functions $f\in C^k_b(X)$ such that 
$D^kf \in C^\alpha(X;\mathcal L^k(X))$, endowed with the norm
\[
\|f\|_{C_b^{k+\alpha}(X)} := \|f\|_{\infty} + \sum_{j=1}^k \sup_{x\in X}\|D^jf(x)\|_{\mathcal L^j(X)} +  [D^kf]_{C^\alpha(X;\mathcal L^k(X))} . 
\]
Let now  $H\subset X$ be  a Hilbert space continuously embedded in $X$, with inner product $\langle \cdot , \cdot\rangle_H$. 
A function $\varphi :X\to Y$ is   $H$-H\"older continuous if there is $\alpha \in (0, 1)$ such that 
$[\varphi ]_{C^{\alpha}_H(X,Y)}:= \sup_{x\in X, \, h\in H\setminus\{0\}} 
\{ \|\varphi (x+h) - \varphi (x)\|_Y/\|h\|_H^{\alpha}\} <+\infty$. $C^{\alpha}_H(X,Y)$ is the space of the functions in 
$C_b(X, Y)$ that are $H$-H\"older continuous with exponent $\alpha$,   with  norm 
$\|\varphi\|_{C^{\alpha}_H(X,Y)}:= \|\varphi\|_{\infty} + [\varphi]_{C^{\alpha}_H(X,Y)}$. 

A function $\varphi :X\to Y$ is
$H$-differentiable at $x\in X$ if there exists $G\in \mathscr L (H, Y)$ such that 
$\|\varphi (x+h) - \varphi (x) -G(h)\|_Y= o(\|h\|_H)$, as $h\to 0$ in $H$. In this case the operator $G$ is unique, and denoted by $D_H\varphi(x)$. Again, if $Y=\R$   there is a unique $y\in H$ such that $G(h) =\langle y, h\rangle_H$ for each $h\in H$. Such $y$ is denoted by $\nabla_H\varphi(x)$. 
If $\varphi $ is  differentiable at $x$ it is also $H$-differentiable at $x$, and 
\[
\frac{\partial \varphi}{\partial h}(x):= Y- \lim_{t\to 0}\frac{\varphi(x+th) -\varphi(x) }{h}=  
\langle \nabla_H\varphi(x), h\rangle _H = D_H\varphi (x)(h) = D\varphi (x)(h), \quad h\in H. 
\]
If $\varphi :X\to \R$ is $H$-differentiable at every point, and in its turn $D_H:X\to H^*$ is $H$-differentiable at 
$x\in X$, we set $D^2_H\varphi(x) := D_H(D_H\varphi )(x) \in \mathcal L^2(H)$ (after identifying $\mathscr L(H, H^*)$ 
with $\mathcal L^2(H)$).  The space $C^1_H(X)$ (resp.  $C^2_H(X)$) consists of the continuous, bounded and (resp. twice) 
$H$-differentiable functions such that $D_H\varphi \in C_b(X, H)$ (resp. $D_H\varphi \in C_b(X, H)$ and 
$D^2_Hf\in C_b(X, \mathcal L^2(H))$). 

\subsection{Semigroups of bounded operators in $C_b(X)$}
\label{sgr}
Let $T(t)$ be a semigroup of bounded operators in $C_b(X)$, such that $\|T(t)\|_{\mathscr L(C_b(X))}\leq Me^{\omega t}$ for some $M>0$, $\omega\in \R$ and for every $t\geq 0$. Assume in addition that the function $(t,x)\mapsto T(t)f(x)$ is continuous in $[0, +\infty)\times X$. 

Since we are going to deal with resolvent and spectrum, it is convenient to extend $T(t)$ to the space $C_b(X;\C)$, setting $T(t)(f+ig) = T(t)f + i T(t)g$ for $f$, $g\in C_b(X)$. 

This allows to define a generator  through its resolvent, 
\begin{equation}
\label{risolvente}
R_{\lambda}f(x) := \int_0^{\infty} e^{-\lambda t}T(t)f(x) dt, \quad  \text{Re}\,\lambda > \omega, \; f\in C_b(X; \C), x\in X. 
\end{equation}
Indeed, in the space $C_b(X, \C)$ the family $\{R_{\lambda}:\; $Re $\lambda >\omega \}$ satisfies the resolvent  identity $R_{\lambda}R_{\mu} = (R_\mu -R_{\lambda})/(\lambda -\mu)$ in the half-plane $\Pi := \{\lambda \in \C:$ Re$\,\lambda>\omega\}$,  since $T(t)$ is a semigroup. Moreover, such identity implies that  if $R_{\mu}f =0$ for some $\mu\in \Pi$ then $R_{\lambda }f =0$ for every $\lambda\in \Pi$. In particular, for every $x\in X$ the Laplace transform $G$ of the function $g(t):= e^{-\omega  t}T(t)f(x)$ vanishes for Re$\,\lambda >0$; since $g\in  C_b([0, +\infty))$ then $g(0) =\lim_{\lambda \to \infty} \lambda G(\lambda) =0$, so that $f(x) = g(0) = 0$. Therefore, $R_{\mu}$ is one to one for every $\mu\in \Pi$, and by e.g.   \cite[\S VIII.4]{Y} there exists a unique closed operator whose resolvent operator is $R_{\mu}$ for every $\mu$ with 
Re$\,\mu >\omega$. 
  The part $L$ of such operator in $C_b(X)$ preserves $C_b(X)$ and it is called {\em generator} of $T(t)$ in $C_b(X)$, although it is not an infinitesimal generator in the classical sense. 

From the definition it follows  $T(t)L=LT(t)$ on $D(L)$. For every $x\in X$, the continuity of $T(\cdot)f(x)$ in $[0, +\infty)$ yields easily its differentiability provided $f\in D(L)$, see e.g. \cite[Prop. 4.2]{Cerrai}. 

We recall that a  Borel probability measure $\mu$ in $X$ is called {\em invariant} for $T(t)$ if
\begin{equation}
\label{invariance}
\int_X T(t)f\,d\mu = \int_X f\,d\mu, \quad t>0 , \; f\in C_b(X). 
\end{equation}

\subsection{Gaussian measures}

We list here notation and results  that will be used in the paper, referring to \cite{Boga} for their proofs and for the general theory. 

A probability measure $\gamma$ 
on $(X,\borel(X))$ is {\em Gaussian} if $\gamma\circ f^{-1}$ (defined as $\gamma\circ f^{-1}(B) = \gamma(f^{-1}(B))$ for every   $B\in \mathcal B( \R)$) is a Gaussian measure in $\R$ for every 
$f\in X^*$. The measure $\gamma$ is called {\em centered}   if all the measures 
$\gamma\circ f^{-1}$ have zero mean,  and it is called {\em nondegenerate} if for any $f\neq 0$ the measure 
$\gamma\circ f^{-1}$ is absolutely continuous with respect to the Lebesgue measure. 

We fix a centered Gaussian measure $\gamma$. By the Fernique Theorem, see \cite[Thm. 2.8.5]{Boga}, $\gamma$ has 
finite moments of any order.  
For every $g\in X^*$ the mapping $R:X^*\to \R$, $Rf:= \int_X f(x)g(x)\ \gamma(dx)$ belongs to $X^{**}$, and even if 
$X$ is not reflexive there exists a unique $y\in X$ such that $Rf = f(y)$, for every $f\in X^*$. We set $y=Qg$. The operator $Q\in \mathscr L(X^*, X)$ is called {\em covariance operator}, it is symmetric and positive. 

If $X$ is a Hilbert space we identify as usual $X$ and $X^*$, and therefore $Q\in \mathscr L(X )$ is defined by 
$$
\langle Qx_0, y_0 \rangle  = \int_X \langle x_0, x\rangle \, \langle y_0, x\rangle\,\gamma(dx), \quad x_0, \;y_0 \in X. 
$$
Moreover, $Q$ belongs to ${\mathscr L}_1(X)$. Conversely, if a linear self-adjoint nonnegative operator $Q$ is 
nuclear, then it is the covariance of a centered Gaussian measure, called $\mathcal N_{0, Q}$. 

Let us go back to general Banach spaces.   The closure   of $X^*$ in $L^2(X, \gamma)$ is denoted by $X^*_{\gamma} $. 
For every $g\in X^*_{\gamma} $, the mapping $R$ defined above still have the representation $Rg = g(y)$, for a suitable (unique) $y\in X$, and we set $y=: R_{\gamma}g$. So, $R_{\gamma}$ is the natural extension of $Q$ to the whole $X^*_{\gamma} $. 

The {\em Cameron-Martin space} $H$ consists of the elements $h\in X$ such that the   measure 
$\gamma_h(B):=\gamma(B-h)$, $B\in \mathcal B(X)$,  is absolutely continuous with respect to $\gamma$.   An important characterization of $H$, that yields a Hilbert space structure in it, is the following: we have 
$H= R_\gamma(X^*_\gamma)$, namely $h\in X$ belongs to $H$ if and only if there is 
$\widehat{h}\in X^*_\gamma$ such that $\int_X\widehat{h}(x) g(x)\,\gamma(dx) = g(h)$ for every $g\in X^*$.
In this case, 
$\|h\|_H=\|\widehat{h}\|_{L^2(X,\gamma)}$. Therefore $R_\gamma:X^*_\gamma\to H$ is an isometry,  and
$H$ is a Hilbert space with the inner product 
$\langle h,k\rangle _H := \langle \widehat h, \widehat k\rangle_{L^2(X,\gamma)}$ 
whenever $h=R_\gamma\widehat h$, $k=R_\gamma \widehat k$.

For every $h\in H$,  the density of $\gamma_h$ with respect to $\gamma$ is given by $e^{-\|h\|^2_H/2 + \widehat h}$. It yields the integration by parts formula
\begin{equation}
\label{parti}
\int_X \frac{\partial \varphi }{\partial h}\, \psi \, \gamma(dx) = - \int_X \varphi  \frac{\partial \psi }{\partial h}\,  \gamma(dx) + \int_X  \varphi\, \psi \,  \widehat h \,\gamma(dx) , \quad \varphi, \;\psi \in C^1_b(X). 
\end{equation}
Moreover for every $h\in H$ the function $\widehat h$ is a real Gaussian random variable with law $\mathcal N_{0, \|h\|_H^2}$. In particular, $\widehat h \in L^q(X, \gamma)$ for every $q\in [1, \infty)$ and 
$\| \widehat h\|_{L^q(X, \gamma)} = ( \int_{\R} |\xi|^{q} \mathcal N_{0,1}(d\xi) )^{1/q} \|h\|_H =: c_q \|h\|_H$. 
 
Recalling that for $f\in X^*$ we have $\int_X f(x) g(x)\,\gamma(dx) = g(Qf)$ for every $g\in X^*$, we see that $H=H_Q$ (the space introduced in Subsection (i)), with the same inner product. 
More precisely, referring to the construction of $H_Q$ in \cite[Chapter III]{VTC} and the operators $A$ involved there, we can take $A:X^*\to X^*_\gamma$, 
$Af = f$, so that $A^*= R_\gamma$.  

If $X$ is a Hilbert space,  the Cameron-Martin space is the range of $Q^{1/2}$, and we have precisely $\langle Q^{1/2}x, Q^{1/2}y\rangle_H = \langle x, y\rangle $ for every $x$, $y\in X$, or
 equivalently $\langle h, k\rangle_H = \langle Q^{-1/2}h, Q^{-1/2}k\rangle $.  
 
 If  $\{e_j:\;j\in \N\}$ is any orthonormal basis of $X$  such that $Qe_j =\lambda_je_j$ for every $j\in \N$, then for every $h\in H$ the function $ \widehat h$ may be represented as 
$ \widehat h(x) = \sum_{j:\,\lambda_j\neq 0}\lambda_j^{-1}\langle h, e_j\rangle \langle x, e_j\rangle $. 
The series converges in $L^p(X, \gamma)$ for every $p\in [1, +\infty)$ 
and it converges pointwise only for $x\in H$, in which case we have $ \widehat h(x)= \langle h, x\rangle_H$. For this reason 
$ \widehat h$ is called $\langle  Q^{-1/2}h, Q^{-1/2}\cdot \rangle $ in \cite{DPZrosso,DPZbrutto}. 

We warn the reader that in the literature about Gaussian measures the expression ``reproducing kernel Hilbert space" is used both for $H$ and for $X^*_\gamma$.  

We denote by $\mathcal{FC}^k_b(X)$ the space of the cylindrical functions $f:X\to \R$ such that $f(x)=\varphi(f_1(x),\ldots,f_n(x))$ with 
 $f_1,\ldots,f_n\in X^*$ and $\varphi\in C^k_b(\R^n)$. Any such functions is $k$ times Fr\'echet differentiable, and we have $Df(x)  = \sum_{j=1}^n D_j\varphi
(f_1(x),\ldots,f_n(x))f_j$, $\nabla_Hf(x) = QDf(x)$. Using \eqref{parti}, one proves that for every $p\in [1, \infty)$ and $k\in \N$, the operator $\nabla_H: \mathcal{F}C^k_b(X)\subset L^p(X, \gamma)\to  L^p(X, \gamma ; H)$ is closable, and the domain  of   its closure (still denoted by $\nabla_H$)  is a Banach space endowed with the graph norm, independent of $k$, called  $W^{1,p}(X, \gamma)$.   
Moreover for $k\geq 2$ the operator $(\nabla_H, D_H^2):  \mathcal{F}C^k_b(X)\subset L^p(X, \gamma)\to  L^p(X, \gamma ; H)\times L^p(X, \gamma ; \mathscr L_2(H))$ is closable too, and the domain of its closure, endowed with the graph norm, is independent of $k$ and called $W^{2,p}(X, \gamma)$. 

The negative of the formal adjoint of $D_H$ is called {\it Gaussian divergence} and denoted by div$_\gamma$. More precisely, a vector field $F\in L^1(X, \gamma;H)$ has Gaussian divergence if there exists (a unique) $\beta\in L^1(X, \gamma)$ such that 
$\int_X\langle \nabla_H \varphi, F\rangle_H\,\gamma(dx) = \int_X \varphi(x) \beta(x)\,\gamma(dx)$, for every $\varphi \in \mathcal{F}C^1_b(X)$. In this case we set  div$_\gamma f:= -\beta$. 
 
\section{The classical O-U semigroup}
\label{sect:OUclassico}
 
Here $X$ is a separable Banach space endowed with a centered Gaussian measure $\gamma$. The proofs of the statements of this section may be found in the book \cite{Boga}, unless otherwise specified. 

The  Ornstein--Uhlenbeck semigroup is  defined through the Mehler formula by
\begin{equation}
\label{OU}
T(0)f = f, \quad T(t) f(x) := \int_X f(e^{-t}x+ \sqrt{1-e^{-2t}}y) \gamma(dy), \quad t>0, \; f\in  C_b(X).
\end{equation}
It is a contraction semigroup in $C_b(X)$,  and $\gamma$ is its unique invariant measure. 
It  is not strongly continuous,  not even in   $BUC(X)$. In fact, it is easily seen that for every $f\in BUC(X)$ we have 
$\lim_{t\to 0^+} \|T(t)f-f\|_{\infty} =0$ if and only if 
$\lim_{t\to 0^+} \|f(e^{-t}\cdot) - f\|_{\infty} =0$. 
 However,  for every $f\in C_b(X)$ the function $(t,x)\mapsto T(t)f(x)$ is continuous in 
$[0, \infty)\times X$ by the Dominated Convergence Theorem, and this allows to define the generator $L$ as in Section \ref{SecPrelimin}(iii). 
Moreover, $T(t)$ is strongly continuous in the
mixed topology, see \cite{GK,GvN}.   

Coming back to the norm topology, $T(t)$ is not analytic and even not continuous in norm in $(0, \infty)$, since 
$\|T(t) - T(s)\|_{\mathscr L(C_b(X))} \geq 2$ for $t\neq s$, as a consequence of \cite[Prop. 2.4]{PvN}. 
The semigroup 
$T(t)$ is   smoothing  along the Cameron-Martin space $H$. More precisely, for every $f\in C_b(X)$ and $t>0$, $T(t)f$ is 
$H$-differentiable at every $x\in X$, and we have
\begin{equation}
\label{funzionegradiente}
\langle \nabla_H T(t)f(x), h\rangle_H =  \frac{e^{-t}}{\sqrt{1-e^{-2t}}} \int_X f(e^{-t}x + \sqrt{1-e^{-2t}}y) \hat{h}(y) \gamma(dy), \quad h\in H. 
\end{equation}
Therefore, using the H\"older inequality and recalling that $\|\hat{h}\|_{L^1(X, \gamma)} \leq \|h\|_H$,   $\|\hat{h}\|_{L^q(X, \gamma)} = c_q\|h\|_H$, 
for every $f\in C_b(X)$ and $x\in X$  we get 
\begin{equation}
\label{stimagradienteH}
\begin{array}{rl}
(i) &\|\nabla_H T(t)f(x)\|_H \leq  \frac{e^{-t}}{\sqrt{1-e^{-2t}}}  \|f\|_{\infty}, 
\\
(ii) &\|\nabla_H T(t)f(x)\|_H \leq  \frac{c_{p'}e^{-t}}{\sqrt{1-e^{-2t}}} [(T(t)|f|^p)(x)]^{1/p},   \quad  p\in (1, \infty), 
\end{array}
\end{equation}
and moreover $ \nabla_H T(t)f :X \to H$ is continuous. 
If in addition  $f\in C^1_b(X)$, then  $T(t)f\in C^1_b(X)$ for any $t\geq 0$, and   
\begin{equation}
\label{derivatedirezT(t)f}
\frac{\partial  T(t)f}{\partial h} (x) = DT(t)f(x)(h) = e^{-t}T(t)(Df(\cdot)(h)),   \quad x, \;h\in X, 
\end{equation}
so that $\sup_{x\in X}\|DT(t)f(x)\|_{X^*} \leq e^{-t}\sup_{x\in X}\|Df(x)\|_{X^*} $. Iterating, we get 
 $T(t) C^k_b(X)\subset  C^k_b(X)$ for any $t\geq 0$, $k\in \N$, and   $\sup_{x\in X}\|D^kT(t)f(x)\|_{\mathcal L^k(X)} \leq e^{-kt}\sup_{x\in X}\|D^kf(x)\|_{\mathcal L^k(X)} $. 
 
Notice that \eqref{funzionegradiente} and  \eqref{stimagradienteH} 
describe a smoothing property of $T(t)$, while the subsequent statements assert that 
$T(t)$ preserves the spaces $C^k_b(X)$ and it is contractive there. However,  $T(t)$ regularizes only along $H$ and it does 
not map $C_b(X)$ into $C^1 (X)$. 

The continuity of $\nabla_HT(t)f$ for $f\in C_b(X)$ and estimate \eqref{stimagradienteH}(i) yield the embedding $D(L)\subset C^1_H(X)$ through the representation formula \eqref{risolvente} for $R(\lambda,L)$. Here, $L$ is the generator of $T(t)$  defined in Section 2(\ref{sgr}). Moreover,  for every $f\in D(L)$, $D_Hf\in C^{\theta}_b(X, H)$ for every $\theta\in (0, 1)$, and it also satisfies a Zygmund condition along $H$, see
\cite{CL}.  A Schauder type theorem was proved in \cite{CL} for $H$-H\"older continuous functions, and precisely: for every $\alpha \in (0, 1)$, $\lambda >0$ and $f\in C^{\alpha}_H(X)$, $R(\lambda, L)f\in C^2_H(X)$ and $D^2_H R(\lambda, L)f\in C^{\alpha}_H(X, \mathcal L^2(H))$. 

The semigroup $T(t)$ is readily extended  to $L^p(X,\gamma)$, for  every   $p\in [1, \infty)$. Indeed,   we have
\begin{equation}
\label{maggLp}
 \int_X |T(t)f(x)|^p\,\gamma (dx) \leq \int_X T(t)(|f |^p)\,\gamma(dx) = \int_X |f |^p \gamma (dx) , \quad t>0, \; f\in C_b(X),  
 \end{equation}
by the H\"older inequality and the invariance of $\gamma$. Hence $\{T(t):\; t\geq 0\}$ is uniquely extendable to a contraction semigroup $\{T_p(t):\; t\geq 0\}$ in $L^p(X,\gamma)$. Moreover,
\begin{itemize}
\item[(i)] $\{T_p(t):\; t\geq 0\}$ is  strongly continuous  in $L^p(X, \gamma)$, for every $p\in [1, \infty)$;
\item[(ii)]  
$T_2(t)$ is self-adjoint and nonnegative in $L^2(X, \gamma)$ for every $t>0$; 
\item[(iii)] $\int_X T_p(t)f\,\gamma(dx) = \int_Xf\,\gamma(dx) $, for every $f\in L^p(X,\gamma)$; 
\item[(iv)] (hypercontractivity) for any $p,q>1$ and $t>0$  such that $q\leq 1+ (p-1)e^{2t}$, $T(t)$ maps $L^p(X, \gamma)$ into $L^q(X, \gamma)$ and 
$\|T(t)f\|_{L^q(X, \gamma)} \leq  \|f\|_{L^p(X, \gamma)}$ for every $f\in L^p(X, \gamma)$. For $q>1+ (p-1)e^{2t}$, $T(t)(L^p(X, \gamma))$ is not contained in $L^q(X, \gamma)$. 
\end{itemize}

For $p\in (1, \infty)$, $L^p$  estimates for $\|\nabla_HT_p(t)f\|_H$ are obtained similarly to \eqref{maggLp}. For every $f\in C_b(X)$,  \eqref{stimagradienteH}(ii) yields 
$$ 
\int_X \|\nabla_HT(t)f(x)\|_H^p\,\gamma (dx) \leq \frac{c_{p'}e^{-t}}{\sqrt{1-e^{-2t}}} \int_X T(t)(|f |^p)\,\gamma(dx) =  \frac{c_{p'}e^{-t}}{\sqrt{1-e^{-2t}}} \int_X |f |^p \gamma (dx) . 
$$
This argument fails for $p=1$, since  \eqref{stimagradienteH}(ii) holds only for $p>1$. Indeed, $T(t)$ does not map $L^1(X, \gamma)$ into $W^{1,1}(X, \gamma)$ for $t>0$, 
even in the simplest case $X=\R$ where $\gamma$ is the standard Gaussian measure (see for instance
\cite[Corollary 5.1]{MPP}). 
For $1\leq p<\infty$, using   formulae \eqref{derivatedirezT(t)f} in $C^1_b(X)$, one obtains that 
 $T_p(t)$ preserves  $W^{1,p}(X, \gamma)$ for every $t>0$, and $\|T_p(t)f\|_{W^{1,p}(X, \gamma)} \leq \| f\|_{W^{1,p}(X, \gamma)} $
for every $f\in W^{1,p}(X, \gamma)$. 

Let us denote by  $L_p$ the infinitesimal generator of $T_p(t)$ in 
$L^p(X,\gamma)$. It is not hard to see that every  $f\in  \mathcal{F}C^2_b(X)$ belongs to $ D(L_p)$, and using \eqref{parti} we get 
\begin{equation}
\label{divgrad}
L_pf(x)={\rm div}_\gamma \nabla_Hf (x) = \sum_{j=1}^{\infty}
\Bigl(\partial^2_{h_jh_j}f(x) - \hat{h}_j(x) \partial_{h_j}f(x)\Bigr), \quad \gamma{\rm -a.e.}\; x\in X, 
\end{equation}
where $\{h_j:\,j\in \N\}$ is any orthonormal basis of $H$. Moreover, $\mathcal{F}C^2_b(X)$ is a core of $L_p$ for every $p\in [1, \infty)$. 
In other words, $D(L_p)$ consists of all $f\in L^p(X, \gamma)$ such that there exists a sequence 
$(f_n)$ in $ \mathcal{F}C^2_b(X)$ which converges to $f$ in $L^p(X, \gamma)$ and such that 
$L_pf_n = {\rm div}{_\gamma} \nabla_H f_n$ converges in $L^p(X, \gamma)$. The Meyer inequalities, see \cite{Meyer}, yield 
\begin{equation}
\label{Meyer}
D(L_p)=W^{2,p}(X,\gamma), \quad 1<p<\infty, 
\end{equation}
with equivalence of the respective norms (an independent analytic proof is in \cite[Section 5.5]{Boga}). For $p=2$,  $L_2$ is the operator associated with the Dirichlet form
\begin{equation}\label{DirForm}
{\cal D}(f,g) = \int_X \langle \nabla_Hf,\nabla_Hg\rangle _H\ d\gamma, \qquad f,g\in W^{1,2}(X,\gamma), 
\end{equation}
namely we have 
\begin{align*}
D(L_2) = \{u\in W^{1,2}(X,\gamma): &\exists\ f\in L^2(X,\gamma)\ \text{ s. t. }\ 
{\cal D}(u,g)= - \langle f,g\rangle_{L^2(X,\gamma)}
\\ 
&\forall\ g\in W^{1,2}(X,\gamma)\},\quad L_2u = f. 
\end{align*}
A complete description of the spectral properties of $L_2$ is available. Even more, there is an explicit orthonormal basis of $L^2(X, \gamma)$ made by eigenfunctions of $L_2$. Such eigenfunctions are the Hermite polynomials, defined for every multiindex $\alpha \in \Lambda:=\{ \alpha \in ( \N \cup \{0\})^{\N} , \;\alpha = (\alpha_j), \; |\alpha |= \sum_{j=1}^{\infty} \alpha_j < \infty\}$, by
 \begin{equation}
\label{Hermiteinf}
H_{\alpha}(x) = \prod_{j=1}^{\infty} H_{\alpha_j}(\hat{h}_j(x)), \quad x\in X, 
\end{equation}
where for $k\in  \N \cup \{0\}$, $H_{k}:\R \to \R$ is the  polynomial  $ H_k(\xi)=\frac{(-1)^k}{\sqrt{k!}}\exp\{\xi^2/2\}\frac{d^k}{d\xi^k}\exp\{-\xi^2/2\}$, for every $\xi\in \R$. 
 
All the polynomials $H_\alpha$ belong to $L^p(X,\gamma)$ for every $p\in [1, \infty)$,  and the set 
$\{ H_{\alpha}: \; \alpha \in \Lambda\}$ is an orthonormal basis of $L^2(X, \gamma)$. 
Moreover, denoting by ${\mathcal X}_k$ the closure of span $\{ H_{\alpha}: \; \alpha \in \Lambda, \; |\alpha| = k\} $ in $L^2(X, \gamma)$, 
we have the {\em Wiener chaos decomposition}, 
\[
L^2(X, \gamma) = \bigoplus _{k\in  \N \cup \{0\}} {\mathcal X}_k. 
\]
The spectrum of 
$L_2$ is equal to $-\N \cup \{0\}$. For every $k\in \N \cup\{0\}$, ${\mathcal X}_k$ is the eigenspace 
of $L_2$ with eigenvalue $-k$.  ${\mathcal X}_0$ is the kernel of $L_2$, consisting of constant functions, and 
${\mathcal X}_1= X^*_\gamma$. 

\section{Ornstein-Uhlenbeck semigroups in Hilbert spaces}
\label{sect:Hilbert}

Here $X$ is a separable real Hilbert space, $Q\in \mathscr L(X)$ is a self-adjoint nonnegative operator, and 
$A:D(A)\subset X\to X$ is the infinitesimal generator of a strongly continuous semigroup $e^{tA}$ in $X$.     

We consider the Ornstein-Uhlenbeck operator formally defined by
\begin{equation}
\label{Lformale}
\mathcal L u(x) = \frac{1}{2} \text{Tr}[QD^2u(x)] + \langle  Ax, \nabla u(x)\rangle  . 
\end{equation}
The standing assumption of  this section  is that the linear operators $Q_t$ defined by  
\begin{equation}
\label{Qt}
Q_tx = \int_0^t e^{sA}Qe^{sA^*}x\,ds, \quad t>0, \;x\in X, 
\end{equation}
are nuclear  ($Q$ itself does not need to be nuclear).
Under such assumption, in \cite[Ch. 6]{DPZbrutto}  it was proved that for very good initial data, namely $f\in BUC^2(X)$ such that $QD^2f \in BUC(X; \mathscr L_1(X))$, the initial value problem 
 \begin{equation}
  \label{P}
 u_t(t,x)  = \mathcal L u(t, \cdot )(x) , \;\; t\geq 0, \;x\in D(A); \quad  u(0,x) = f(x), \;\;x\in X,
  \end{equation}
has a unique strict solution, which is a continuous function $u:[0, +\infty)\times X\to \R$ such that   $u(t, \cdot)\in BUC^2(X)$  for every $t\geq 0$,   $QD^2u(t,x)\in \mathscr L_1(X)$ for every $t\geq 0$ and $x\in X$, $u(\cdot, x) $ is continuously differentiable in $[0, +\infty)$ for every $x\in D(A)$, and satisfies  \eqref{P}. 
Moreover, the solution is given by 
\begin{equation}
\label{soluzione}u(t,x) = \int_X f(e^{tA}x +y) \mathcal \mu_t(dy), \quad t\geq 0, \; x\in X, 
 \end{equation}
where $\mu_t$ is the centered Gaussian measure $ \mathcal N_{0,Q_t}$ with covariance  $Q_t$ for $t>0$, and $\mu_0 =\delta_0$. 

\subsection{Ornstein-Uhlenbeck semigroups in spaces of continuous functions}

The right hand side of \eqref{soluzione} is meaningful for every $f\in {\mathcal B}_b(X)$. Setting
\begin{equation}
\label{OUHilbert}
T(t)f (x) := \int_X f(e^{tA}x +y) \mathcal \mu_t(dy), \quad t\geq 0, \; f\in {\mathcal B}_b(X), \;x\in X, 
 \end{equation}
$T(t)$ is a contraction semigroup in ${\mathcal B}_b(X)$. The fact that $T(t)$ maps ${\mathcal B}_b(X)$ into itself and it is a contraction is obvious. The fact that $T(t)$ is a semigroup is less obvious. It can be proved rewriting $T(t+s)$, for $t$, $s>0$, as 
$$
T(t+s)f(x) =  \int_X f(e^{(t+s)A}x +w)  (\mu_t \circ (e^{sA})^{-1}\star \mu_s)(dw), \quad f\in {\mathcal B}_b(X), \;x\in X, 
$$
and checking that $\mu_t \circ (e^{sA})^{-1}\star \mu_s = \mu_{t+s}$, or else recalling that $T(t)$ is the transition semigroup of the stochastic differential equation 
\begin{equation} 
\label{stoch}
dX_t = AX_t\,dt + \sqrt{Q}\,dW_t, \;\;t>0, \quad X(0) = x, 
\end{equation}
where $W_t$ is any cylindrical Wiener process in $X$. Indeed, for every $x\in X$  the unique mild solution of \eqref{stoch} is $X_t= e^{tA}x +\int_0^te^{(t-s)A}Q^{1/2}dW_s$, and the law of the stochastic convolution $\int_0^te^{(t-s)A}Q^{1/2}dW_s$ is precisely $\mathcal N_{0,Q_t}$. See \cite[Ch. 5]{DPZrosso}. Therefore, 
$$
T(t)f(x) = {\mathbb E}(f(X_t)), \quad t\geq 0, \;f\in {\mathcal B}_b(X), \; x\in X. 
$$
If $A=-I$ and $Q$ is nuclear, setting  $\gamma := \mathcal N_{0, 2Q}$, $T(t)$ coincides with the classical Ornstein-Uhlenbeck semigroup considered in  Section \ref{sect:OUclassico}.  If $A=0$, $T(t)$ may be called {\em heat semigroup}. In this case, $Q_t =tQ$ so that setting $y=\sqrt{t}z$ in the right hand side of \eqref{OUHilbert} we get a simpler representation formula for $T(t)$, 
$$
T(t)f (x) := \int_X f(x +\sqrt{t} z)  \mu(dz), \quad  t\geq 0, \; f\in {\mathcal B}_b(X), \;x\in X, 
$$
where $\mu :=  \mathcal N_{0, Q}$. 

Going back to the general case, the representation formula \eqref{OUHilbert} yields that $T(t)$ is a Feller semigroup, namely it maps $C_b(X)$ into itself  and, in fact, it maps the subspaces $BUC(X)$, $C^{\alpha}_b(X)$, $C^{ k}_b(X)$, $C^{\alpha +k}_b(X)$ into themselves, for $\alpha\in (0, 1)$, $k\in \N$. 
In particular, for $f\in C^{1}_b(X)$ we have  
\begin{equation}
\label{comm}
\langle \nabla T(t)f(x), h\rangle = \int_X \langle e^{tA^*}\nabla f(e^{tA}x+y), h\rangle \,\mu_t(dy), \quad x, \;h\in X. 
\end{equation}
$T(t)$ is strong-Feller (namely, it maps ${\mathcal B}_b(X)$ into $C_b(X)$)  iff (see also \cite[Remark 9.20]{DPZrosso})
\begin{equation}
\label{controllability}
e^{tA}(X) \subset Q_t^{1/2}(X), \quad t>0. 
\end{equation}
In this case, $T(t)$ maps ${\mathcal B}_b(X)$ into $C^k_b(X)$ for every $k\in \N$ (\cite[Ch. 6]{DPZbrutto}), and  the operators
\begin{equation}
\label{Lambda_t}
\Lambda_t = Q_t^{-1/2}e^{tA}, \quad t>0, 
\end{equation}
play an important role in the rest of the theory. 
First, $\Lambda_t \in \mathscr L(X)$ for every $t>0$. 
Moreover, 
for every $k\in \N$ there exists $C_k>0$ such that 
\begin{equation}
\label{regolarizzaz}
\| D^kT(t)f(x)\|_{\mathcal L^k(X)} \leq C_k\|\Lambda_t\|_{\mathscr L(X)}^k\|f\|_{\infty},  \quad t>0, \; 
f\in {\mathcal B}_b(X), \; x\in X. 
\end{equation}
A proof for $k=1, 2$ is in \cite[Ch. 6]{DPZbrutto}. For general $k$, \eqref{regolarizzaz} follows e.g. from \cite[Sect. 5.1, Prop. 3.3(ii)]{LR}. 
 
Notice that \eqref{controllability} is not satisfied if $A=-I$, and, more generally, if $A$ generates a strongly continuous group. Instead, it is satisfied if $Q=I$, and in this case $\|\Lambda_t\|_{\mathscr L(X)} \leq Me^{\omega t}t^{-1/2}$ for some $M>0$, $\omega\in \R$, and for every $t>0$. See  \cite[Appendix B]{DPZbrutto}, \cite[Thm. 3.5(3)]{GvN}. 

Anyway, smoothing properties along  $H:=Q^{1/2}(X)$ are available also in the case where $H$ is properly contained 
in $X$, provided that $e^{tA}$ maps $H$ into itself, and that $S_H(t):= e^{tA}_{|H}:H\to H$ is a strongly continuous 
semigroup in $H$. In this case $e^{tA}$ maps $H$ into $Q_t^{1/2}(X)$ for every $t>0$, and 
$\sup_{0<t<1}\|e^{tA}\|_{\mathscr L(H, Q_t^{1/2}(X))} <\infty$, by \cite[Thm. 3.5]{GvN}. This allows to prove that $T(t)$ 
is smoothing along $H$, by arguments similar to the ones that led to \eqref{stimagradienteH}(ii). 
See  \cite[Sect. 2]{MaasVN4}, and  \cite{LR} for representation formulae and estimates for any order $H$-derivatives of 
$T(t)f$ when $f\in C_b(X)$. 

Let us discuss strong continuity. Even in the case $X=\R$,  $T(t)$ is not strongly continuous in $BUC(X)$ unless $A=0$ 
(let alone in $C_b(X)$).  However, 
it is not hard to show that $\mu_t$ converges weakly to $\delta_0$ as $t\to 0$ (namely, $\lim_{t\to 0}\int_X f(y)\mu_t(dy) = f(0)$ for every $f\in C_b(X)$) and this implies 
$$
\lim_{t\to 0}\|T(t)f - f(e^{tA}\cdot)\|_{\infty} =0, \quad f\in BUC(X). 
$$
So, the subspace $BUC_S(X)$ of strong continuity of $T(t) $ in $BUC(X)$ is $\{f\in BUC(X):\; \|T(t)f - f(e^{tA}\cdot)\|_{\infty} \to 0$ as $t\to 0\}$. If \eqref{regolarizzaz} holds, $T(t)(C_b(X) )\subset   BUC(X)$ for every $t>0$ and therefore $BUC_S(X)$ coincides with the subspace of strong continuity of $T(t) $ in $C_b(X)$. In the general case, the subspace of strong continuity of $T(t) $ in $C_b(X)$ is not known. 
However, $T(t)$ is strongly continuous in $C_b(X)$ with respect to the mixed topology, see \cite{Cerrai,GK}. In particular, the function $(t,x)\mapsto T(t)f(x)$ is continuous in $[0, +\infty) \times X$ for every $f\in C_b(X)$, and this allows to define a generator $L$ as in Section   \ref{SecPrelimin}(iii). 
Moreover, setting $\Delta_hf = (T(h)f-f)/h$ for $h>0$,  we have
\begin{align} \nonumber
D(L) &= \{f\in C_b(X):\; \exists g\in C_b(X) \; {\rm s.t.} \lim_{h\to 0} \|\Delta_hf -g\|_{\infty} =0, \; 
\limsup_{h\to 0} \|\Delta_hf \|_{\infty}<+\infty\}, 
\\ \label{dominio_generatore}
&= \{f\in C_b(X):\; \exists g\in C_b(X) \; {\rm s.t.} \lim_{h\to 0} \Delta_hf = g \;\text{ in the weak topology} \}, 
\quad  Lf  = g.
\end{align}
See \cite{GK,GvN}. An analogous characterization with the space $C_b(X)$ replaced by $BUC(X)$ is in \cite{Priola}. Still in  \cite{GK} it was proved that (similarly to the case of strongly continuous semigroups in Banach spaces) any subspace $D\subset D(L)$ which is dense in $C_b(X)$ in the mixed topology and such that $T(t)(D)\subset D$, is a core for $L$, namely for every $f\in D(L)$ there exists a net $(f_\alpha)\subset D(L)$ such that $f_\alpha\to f$ and $Lf_\alpha \to Lf$ in the mixed topology. Using this property it is possible to show that 
\begin{align}
\mathcal{F }_0 := \{  f\in C_b(X):\: &f = \varphi(\langle \cdot, a_1\rangle, \ldots, \langle \cdot, a_n\rangle); \; \varphi \in C^{2}_b(\R^n), 
\nonumber \\ \label{F0}
&\;n\in \N, \; a_i\in D(A^*),\ \langle \cdot , A^*\nabla f \rangle\in C_b(X)\}
\end{align}
and its subspace
 $\mathcal{F}C^{\infty}$ (whose members are the functions $f$ represented as in \eqref{F0} with 
 $\varphi \in C^{\infty}_c(\R^n)$)  
are cores for $L$; moreover    we have (\cite[Thm. 4.5]{GK}, \cite[Thm. 6.6]{GvN})  
\begin{equation}
\label{L}
Lf (x) =  \frac{1}{2} \text{Tr}[QD^2f(x)] + \langle  x, A^*\nabla f(x)\rangle , \quad  f\in \mathcal{F }_0, \;x\in X, 
\end{equation}
where the right-hand side is equal to $\mathcal Lf(x)$ for every $x\in D(A)$. 
Related results with $BUC(X)$ replacing $C_b(X)$ are in \cite{Cerrai,CerraiGozzi,Priola}. 
In some papers, see e.g. \cite{GK}, also the realization of $T(t)$ in the {\em weighted} spaces 
$C_m(X)= \{ f\in C(X; \R):\; \|f\|_{C_m(X)} :=\sup_{x\in X}|f(x)|/(1+\|x\|^m)<\infty\}$ has been studied. 
 
In finite dimension $T(t)$ is analytic iff $A=0$. Instead, if $X$ is infinite dimensional, we have 
$\|T(t)-T(s)\|_{{\mathscr L } (BUC_S(X))} =2$ and therefore $\|T(t)-T(s)\|_{{\mathscr L } (C_b(X)} \geq 2$ 
whenever $\mu_t$ and $\mu_s$ are singular (which is the case for every $t$, $s>0$ if $A=0$). The same equality holds if $e^{tA}\neq e^{sA}$, see 
\cite{PvN,VanNervZabc}. Therefore, $T(t)$ is not norm continuous, and hence not analytic, both in the case $A=0$ and  in the case $A\neq 0$. See \cite{GvN,PvN,VanNervZabc}. 

An alternative proof of norm discontinuity in the case $A=0$ comes from \cite{MRS}, where it has been proved that the spectrum of  the part of $L$  in $BUC(X;\C)$ is the halfplane $\{ \lambda\in\C:$ Re$\,\lambda \leq 0\}$, and for every $t>0$ the spectrum of $T(t)$ in $BUC(X;\C)$ is the whole closed unit disk. 

Schauder type results in the usual H\"older spaces are available if \eqref{controllability} holds, under the further assumption
\begin{equation}
\label{holder}
\exists \; M, \;\theta >0, \;\omega \in \R:\quad   \|\Lambda_t \|_{\mathscr L(X)} \leq \frac{Me^{\omega t}}{t^{\theta}}, \quad t>0. 
\end{equation}
Easy examples such that \eqref{Qt} and \eqref{holder} hold (with any $\theta >0$) are given in \cite[Ex. 6.2.11]{DPZbrutto}. 
The following theorem is taken from \cite[Sect. 5.1]{LR}. 

\begin{theorem}
\label{SchauderLR}
Let \eqref{Qt} and \eqref{holder}  hold. For every $f\in C_b(X)$ and $\lambda >0$, let $u=R(\lambda, L)f$. Then
\begin{itemize}
\item[(i)] If $1/\theta \notin \N$, then $u
  \in C^{1/\theta}_b(X)$, and there is $C>0$, independent of $f$, such that 
 $\|u\|_{C^{1/\theta}_b(X)} \leq C\|f\|_{\infty}$. 
\item[(ii)] If in addition $f\in C^{\alpha}_b(X)$ with $\alpha \in (0, 1)$ and $\alpha +1/\theta  \notin \N$, then $u\in C^{\alpha + 1/\theta }_b(X)$  and there is $C>0$, independent of $f$, such that 
$\|u\|_{C^{\alpha + 1/\theta }_b(X)} \leq C\|f\|_{ C^{\alpha}_b(X)}$. 
\end{itemize}
\end{theorem}
Statement (ii) was already proved in \cite{CerraiLibro} in the case that $A$ is the realization of a second order elliptic system with general boundary conditions in $X =L^2(\Omega)$, $\Omega$ being a bounded open set in $\R^n$, and suitable assumptions on $Q$ that yield $\theta =1/2$. See also  \cite{CDP} for an earlier result. 

Statement (i) implies that $D(L)\subset C^{1/\theta}_b(X)$ if $1/\theta \notin \N$, with continuous embedding. Statement (ii) implies that the domain of the part of $L$ in 
$ C^{\alpha}_b(X)$ is continuously embedded in $C^{\alpha + 1/\theta }_b(X)$ if $\alpha +1/\theta  \notin \N$. In both cases, we gain ``$1/\theta$ degrees" of regularity.  

Both for $\alpha =0$ and for $\alpha >0$, in the critical cases $\alpha +1/\theta = k\in \N$ we cannot expect that $u\in C^k(X)$; in \cite{LR} it is proved that $u$ belongs to a suitable Zygmund space, which is continuously embedded in all  spaces $C^{k-\varepsilon}_b(X)$ with $\varepsilon \in (0, 1)$.   This difficulty arises even in finite dimension, for instance 
if $X=\R^n$, $A=0$, $Q=2I$ we have  $\mathcal L =\Delta $, $Q_t = 2tI$ and  \eqref{holder}  holds with $\theta =1/2$, but if $\lambda u -\Delta u =f\in C_b(\R^n)$ with $n\geq 2$, $u$ is not necessarily a $C^2$ function. 

If  $e^{tA}$ maps $H=Q^{1/2}(X)$ into itself, and $S_H(t)= e^{tA}_{|H}:H\to H$  is a strongly continuous semigroup in $H$,  Schauder  theorems similar to the ones stated in Sect. \ref{sect:OUclassico} were proved in  \cite{LR}:  for every $\alpha \in (0, 1)$, $\lambda >0$ and $f\in C^{\alpha}_H(X)$, $R(\lambda, L)f\in C^2_H(X)$ and $D^2_H R(\lambda, L)f\in C^{\alpha}_H(X, \mathcal L^2(H))$. 

Schauder type regularity results are available also for evolution equations with bounded and continuous data, see \cite{LR}. 
 
The asymptotic behavior of $T(t)$ is well understood if 
\begin{equation}
\label{invmeas}
\sup_{t>0} {\rm Tr}\,(Q_t)  =  \int_0^{\infty} {\rm Tr}(e^{sA}Qe^{sA^*})ds <+\infty . 
\end{equation}
Next statements are taken from \cite[Sect. 11.3]{DPZrosso}, \cite[Sect. 10.1]{DPZbrutto}. If \eqref{invmeas} holds  there exists a nuclear self-adjoint operator $Q_{\infty}$, given by
\begin{equation}
\label{Qinfty}Q_{\infty}x =  \int_0^{\infty}  e^{sA}Qe^{sA^*}x\,ds, \quad x\in X, 
\end{equation}
which maps $D(A^*)$ into $D(A)$  and satisfies the identity (called Lyapunov equation)
\begin{equation}
\label{Lyapunov}
Q_{\infty}A^*x + AQ_{\infty}x= -Qx, \quad x\in D(A^*), 
\end{equation}
Such identity is easily obtained recalling that  $\langle Q_{\infty}e^{tA^*}x, e^{tA^*}y\rangle = \langle Q_{\infty}x, y\rangle -  \langle Q_{t}x, y\rangle$ for every  $x$, $y\in X$. Indeed, taking $x$, $y \in  D(A^*)$, differentiating in time and taking  $t=0$ we get $\langle Q_{\infty}A^*x,  y\rangle + \langle Q_{\infty}x, A^*y\rangle  =\langle Q x, y\rangle$ and \eqref{Lyapunov} follows by the density of $D(A^*)$. 

Moreover,  the  Gaussian measure $\mu_{\infty}:= \mathcal N_{0,  Q_{\infty}}$ is invariant for $T(t)$, namely
$$
\int_X T(t)f(x) \,\mu_{\infty}(dx) = \int_X f(x) \,\mu_{\infty}(dx), \quad t>0, \; f\in C_b(X). 
$$
In fact,  it is possible to show that \eqref{invmeas} holds iff there exists a probability  invariant measure for $T(t)$ iff there exists a self-adjoint nonnegative nuclear operator $P$ mapping $D(A^*)$ into $D(A)$ and such that $  PA^*x  + APx = -Qx$ for every $x\in  D(A^*)$ (which is equivalent to 
 $2 \langle  PA^*x, x\rangle + \langle Qx, x\rangle = 0$ for every $x\in D(A^*)$). Moreover, any  invariant measure is given by $\nu\star \mu_{\infty}$, $\nu$ being a probability invariant measure for the semigroup $R(t)$ defined by $R(t)f(x) = f(e^{tA}x)$ (e.g. \cite{Zab1981},  \cite[Thm. 11.17]{DPZrosso}). 
So, if $R(t)$ has no invariant measures except $\delta_0$, $\mu_{\infty}$ is the unique invariant measure for $T(t)$. In particular, this happens if $\lim_{t\to \infty} e^{tA}x =0$ for every $x$. 

If $ \|e^{tA}\|_{\mathscr L (X)} $ vanishes as $t\to\infty$, namely if there are $M$, $\omega >0$ such that  
\begin{equation}
\label{sgrdecay}
  \|e^{tA}\|_{\mathscr L (X)} \leq Me^{-\omega t}, \quad t>0, 
\end{equation}
it is not hard to see that \eqref{invmeas} holds (e.g., \cite[Thm. 11.20]{DPZrosso}), and therefore $\mu_{\infty}$ is well defined and it is the unique invariant measure for $T(t)$. 

Notice that if $Q$ commutes with $e^{tA}$ for every $t$ and \eqref{sgrdecay} holds, then $A$ is invertible, and if in addition it is self-adjoint then  $Q_{\infty} = -QA^{-1}/2 = -A^{-1}Q/2$. The equality $Q_{\infty} =  -A^{-1}Q/2$ holds even in a more general situation, see the remarks after Theorem \ref{thm:CMG1}. 

It is interesting to compare kernels and  ranges of $Q ^{1/2}$, $Q_t^{1/2}$ and  $Q_{\infty}^{1/2}$ for $t>0$, that  play an important role in the theory. We set
$$
H:=  Q^{1/2}(X), \quad   H_{t}:=  Q_{t}^{1/2}(X), \quad H_{\infty}:=  Q_{\infty}^{1/2}(X), 
$$
endowing them with their natural inner products, described in Sect. \ref{SecPrelimin}(i).  
Using the Lyapunov equation \eqref{Lyapunov} one gets easily  (e.g.,  \cite[Lemma 2.1]{Fuhrman})
$$
e^{tA}H_{\infty}\subset H_{\infty}, \quad \|Q_{\infty}^{-1/2} e^{tA}Q_{\infty}^{1/2}\|_{\mathscr L(X)} \leq 1, \quad t>0. 
$$ 
Therefore,  $e^{tA}_{|H_{\infty}}:  H_{\infty}\to H_{\infty}$ is a contraction semigroup, called $S_{\infty}(t)$. 
Its infinitesimal generator is the part $A_{\infty}$ of $A$ in $H_{\infty}$. 
Since $\langle Q_tx, x\rangle \leq \langle Q_{\infty}x, x\rangle$ for every $t>0$ and $x\in X$, then Ker$\,Q_{\infty} =$ Ker$\,Q_{\infty}^{1/2}\subset \,$Ker$\,Q_t^{1/2}\subset \,$Ker$\,Q ^{1/2} = $ Ker$\,Q $, and 
 $H_t\subset  H_{\infty}$  (we recall that, given self-adjoint operators $T_1$, $T_2\in \mathscr L(X)$, we have $T_1(X)\subset T_2(X)$ iff there exists $C>0$ such that $\|T_1x\|\leq C\|T_2x\|$ for every $x\in X$). 
Instead, the converse inclusion  $H_{\infty} \subset H_t$ is not necessarily satisfied, and by \cite[Prop. 4.1]{Fuhrman} or \cite[Lemma 4]{Ch-MiGol2} it is equivalent to 
\begin{equation}
\label{eq:Fuhrman}
 \|Q_{\infty}^{-1/2} e^{tA}Q_{\infty}^{1/2}\|_{\mathscr L(X)}<1, 
 \end{equation}
namely, to $\|S_{\infty}(t)\|_{{\mathscr L}(H_{\infty})} <1$. 

In the proof of Theorem 11.22 of \cite{DPZrosso}  it was shown that  if  \eqref{controllability} holds, then $H_{\infty}  \subset H_t $ (so that \eqref{eq:Fuhrman} holds) and moreover the operators $Q_t^{-1/2}Q_{\infty}Q_t^{-1/2} -I$ are Hilbert-Schmidt on  $H_{\infty} $, for every $t>0$. Therefore, by the Feldman-Hayek Theorem (e.g. \cite[Thm. 2.23]{DPZrosso}) $\mu_t$ and 
$\mu_{\infty}$ are equivalent measures, for every $t>0$. 

If \eqref{invmeas} holds,  we have (\cite[Thm. 11.20]{DPZrosso})
\begin{equation}
\label{asymptotics}
\lim_{t\to \infty} T(t)f(x) = \int_X f(y)\,\mu_{\infty}(dy), \quad f\in C_b(X), \; x\in X. 
\end{equation}
We notice that if $A=0$, then Tr$\,Q_t = t\,$Tr $Q$, so that \eqref{invmeas} does not hold, and the heat semigroup has no invariant measure. Instead, if $A= -\omega I$ with $\omega>0$, then Tr  $Q_t = (1-e^{-2\omega t} )\,$Tr$\,Q/(2\omega)$, so that  \eqref{invmeas} holds with $Q_{\infty}= Q/(2\omega)$. In particular, as we already mentioned in Section \ref{sect:OUclassico}, the classical Ornstein-Uhlenbeck semigroup has $\gamma $ itself as unique invariant measure (we recall that the covariance of $\gamma$ is $2Q$). 
 
\subsection{Ornstein-Uhlenbeck semigroups in $L^p$ spaces with respect to invariant measures}

Throughout this section we assume that \eqref{invmeas} holds, and we consider $L^p$ spaces with respect to the invariant measure $\mu_{\infty}$, $1\leq p<\infty$. 

For every $f\in C_b(X)$ and $t>0$, the H\"older inequality and the invariance of  $\mu_{\infty}$  yield
$$  
\int_X |T(t)f(x)|^p\mu_{\infty}(dx) \leq \int_X T(t)(|f |^p)\mu_{\infty}(dx) = \int_X |f |^p \mu_{\infty}(dx) 
$$
and therefore, since $C_b(X)$ is dense in $L^p(X, \mu_{\infty})$,  $T(t)$ has a bounded extension to  $L^p(X, \mu_{\infty})$, denoted by $T_p(t)$. The above inequality implies that $T_p(t)$ is a contraction semigroup in $L^p(X, \mu_\infty)$. By the Dominated Convergence Theorem, $\lim_{t\to 0} \|T(t)f -f\|_{L^p(X, \mu_{\infty})} = 0$ for every $f\in C_b(X)$, and this yields 
$\lim_{t\to 0} \|T_p(t)f -f\|_{L^p(X, \mu_{\infty})} = 0$ for every $f\in L^p(X, \mu_{\infty})$.  

The  generator of $T_p(t)$ is denoted by $L_p$. Since $T_p(t)f=T_q(t)f$ for $p\leq q$ and $f\in L^q(X, \mu_{\infty})$, then $L_q$ is the part of $L_p$ in $L^q(X, \mu_{\infty})$, and the subindex $p$ will be written only if needed. 

Notice that, for every $f\in D(L_p)$,  letting $t\to 0$ in the equality $\int_X [(T(t)f -f)/t ]d\mu_{\infty} =0$ we obtain
$\int_X L_pf\, d\mu_{\infty} =0$ . 

Concerning asymptotic behavior, for every $f\in L^p(X, \mu_{\infty})$ we have 
\begin{equation}
\label{asymptoticsLp} 
\lim_{t\to\infty} \left\| T_p(t)f -  \int_X f(y)\,\mu_{\infty}(dy)\right\|_{L^p(X, \mu_{\infty})} =0. 
\end{equation}
If $f\in C_b(X)$, \eqref{asymptoticsLp} is a consequence of  \eqref{asymptotics} through the Dominated Convergence Theorem; if $f\in L^p(X, \mu_{\infty})$ it follows approximating $f$ by a sequence of continuous and bounded functions. 
 
Using the Dominated Convergence Theorem, it is easy to see that the space $\mathcal{F }_0 $ defined in \eqref{F0} is contained in $D(L_p)$ for every $p\in [1, \infty)$, and it is a core for $L_p$ since it is invariant under $T(t)$ and  dense in $L^p(X, \mu_{\infty})$. Another convenient core, used in \cite{DPZbrutto}, is the subspace of $\mathcal{F }_0 $ defined by 
$$
{\mathcal E}_A(X):= \text{span}\; \{\cos (\langle \cdot, h\rangle), \; \sin  (\langle \cdot, k\rangle); \; h, \,k \in D(A^*)\}. 
$$
Necessary and sufficient conditions for $T_2(t)$ be self-adjoint for every $t>0$ (or, equivalently, for  $L_2$ be 
self-adjoint) were given in \cite{Ch-MiGol1} under the assumption that $Q_{\infty}$ is one to one, that was later removed in \cite{GvN}. In both papers, the key tool was the representation of $T_2(t)$ as the second quantization operator of the operator $S_{\infty}(t)^*$,  that goes back to \cite{Ch-MiGol2}.  

\begin{theorem}
\label{thm:CMG1}
The following conditions are equivalent.
\begin{itemize}
\item[(i)] $T_2(t) = T_2(t)^*$ for every $t>0$; 
\item[(ii)] $Q(D(A^*)) \subset D(A)$, and $AQx = QA^*x$ for every $x\in D(A^*)$; 
\item[(iii)] $e^{tA}Q = Qe^{tA^*}$, for every $t>0$; 
\item[(iv)]  $e^{tA}Q_{\infty} = Q_{\infty}e^{tA^*}$, for every $t>0$; 
\item[(v)] $e^{tA}(H)\subset H$, and $S_H(t):= e^{tA}_{|H}:H\to H$ is a self-adjoint strongly continuous semigroup in $H$. 
\end{itemize}
\end{theorem}
We refer to the conditions of Theorem  \ref{thm:CMG1} as ``the symmetric case". 
In such a  case, by the general theory of semigroups  the infinitesimal generator $L_2$ of $T_2(t)$ is self-adjoint too. Moreover $T_2(t)$ is a symmetric Markov semigroup in $L^2(X, \mu_{\infty})$, according to the terminology of \cite{Davies}, and therefore 
 $T_p(t)$ is an analytic semigroup in $L^p(X, \mu_{\infty})$ for every $p\in (1, +\infty)$ with angle of analyticity $\geq \pi(1-|2/p-1|)/2$, by \cite[Thm. 1.4.2]{Davies}. 
In addition, (iv) yields that $Q_{\infty}$ maps $D(A^*)$ into $D(A)$, and on $D(A^*)$ we have $AQ_{\infty} = Q_{\infty}A^*$ ($ = -Q/2$ by the Lyapunov equation). In particular, if $0$ belongs to the resolvent set $\rho(A)$ we get an explicit formula 
for $Q_{\infty} = -\frac{1}{2}A^{-1}Q = -\frac{1}{2}Q(A^*)^{-1}$. 
About condition (v), we remark that $S_H(t)$ is self-adjoint and strongly continuous in $H$ iff $Q^{-1/2}e^{tA}Q^{1/2}$ is self-adjoint and strongly continuous in $X$. Moreover, in the symmetric case not only $S_H(t)$ is strongly continuous, but there are $M_1$, $\beta>0$ such that 
\begin{equation}
\label{decaySH}
\|S_H(t)\|_{\mathscr L(H)} \leq M_1 e^{-\beta t}, \quad t>0. 
\end{equation}
See \cite[Thm. 4.5]{GvN}. Such estimate plays an important role in the asymptotic behavior of $T_p(t)$. 

In the nonsymmetric case, $T_p(t)$ is not in general analytic, even in finite dimension: see the counterexamples in \cite{Fuh1995}.  
Necessary and sufficient conditions for analyticity were studied in the  papers   \cite{Schmuland,Fuh1995,GoldysLincei,GvN,MaasVN1,VanNervZabc}. In particular, \cite{GvN}   contains extensions and improvements of the previous ones, that are summarized in the next theorem. 

\begin{theorem}
\label{Thm:analyticity}
The following conditions are equivalent:
\begin{itemize}
\item[(i)] $T_2(t)$ is an analytic semigroup in $L^2(X, \mu_{\infty})$; 
\item[(ii)] there exists $M>0$ such that  
 $|\langle  Q_{\infty}A^*x, y\rangle | \leq M \langle Qx, x\rangle^{1/2} \langle Qy, y\rangle^{1/2}$, for every $x$, $y\in D(A^*)$;
\item[(iii)] $S_{\infty}(t)$ is an analytic contraction semigroup$^($\footnote{An analytic semigroup $T(t)$ in  a real 
Banach space $\mathcal X$ is called ``analytic contraction semigroup" if there exists a sector 
$\Sigma := \{z\neq 0:\; |$arg$\,z|< \theta\}$ with $\theta >0$ such that the analytic extension $T(z)$ satisfies 
$\|T(z)\|_{\mathcal L(\mathcal X^{\C})}\leq 1$ for every $z\in \Sigma$. $\mathcal X^{\C}$ is the complexification 
of $X$.}$^)$ in $ H_{\infty}$. 
\end{itemize}
If in addition $Q$ has a bounded inverse, the above conditions are also equivalent to 
\begin{itemize}
\item[(iv)] The operator $AQ_{\infty}$ has an extension belonging to $\mathscr L(X)$; 
\item[(v)] The operator $Q_{\infty}A^*$ has an extension belonging to $\mathscr L(X)$. 
\end{itemize}
\end{theorem} 
We refer to the conditions of Theorem  \ref{Thm:analyticity} as ``the analytic case".  

As in the symmetric case, if $T_2(t)$ is  analytic  in $L^2(X, \mu_{\infty})$ then $T_p(t)$ is analytic  in $L^p(X, \mu_{\infty})$ for every $p\in (1,\infty)$, by a simple application of the Stein Interpolation Theorem (e.g., \cite[Sect. 6.2]{Interp}). Moreover, $T_p(t)$ is an analytic contraction semigroup and the optimal angle of analyticity $\theta_p$ has been determined in \cite{MaasVN1}; in \cite{CarDra} it has been proved that such angle coincides with the optimal angle for the bounded $H^{\infty}$ calculus of $-L_p$. 
In addition, in the analytic case the semigroup $e^{tA}$ maps $H$ into itself, and the semigroup 
 $S_H(t)$ is a strongly continuous, bounded analytic semigroup in $H$, see \cite[Thm. 3.3]{MaasVN4}. 
For $p=1$, $T_1(t)$ is not analytic, even in finite dimension. 
Characterizations of the domains $D(L_p)$ as suitable Sobolev spaces are known only in the analytic case.  

The definition of the proper Sobolev spaces relies on the closability of the operator  $\nabla_H: {\mathcal F}_0 \subset L^p(X, \mu_{\infty}) \to L^p(X, \mu_{\infty}; H)$, with $p\in [1, \infty)$. If $f\in  {\mathcal F}_0$, $f(x) = \varphi (\langle x, x_1\rangle, \ldots, \langle x, x_n\rangle)$ with $\varphi\in C^2_b(\R^n)$ and $x_k\in D(A^*)$, we have 
$\nabla_Hf(x) = \sum_{k=1}^n D_k\varphi (\langle x, x_1\rangle, \ldots, \langle x, x_n\rangle) Qx_k$. 
Recalling \eqref{L}, \eqref{parti}  and using the Lyapunov equation  it is easy to see that for  $f$, $g\in {\mathcal F}_0$ we have
\begin{equation}
\label{integrazioneL}
\int_X (Lf \,g  + f \,Lg )\mu_{\infty}(dx) = -\int_X \langle Q\nabla f , \nabla g \rangle \mu_{\infty}(dx) = - \int_X \langle \nabla_Hf , \nabla _Hg \rangle_H \mu_{\infty}(dx). 
 \end{equation}
According to  \cite[Sect. 6]{GGVN}, a sufficient condition for $\nabla_H$ be closable is that $Q$ is one to one and the operator $W: H_{\infty} \to X$, $W(x) = Q^{1/2}Q_{\infty}^{-1/2}x$, is closable in $X$. Another sufficient condition, see 
\cite[Cor. 5.6]{GGVN}, is that $e^{tA}$ maps $H$ into itself  and 
 $S_H(t)$ is strongly continuous in $H$. 
So, in the analytic case (and, in particular, in the symmetric case) $\nabla_H$ is closable in $L^p(X, \mu_{\infty})$ for every $p\in [1, \infty)$. See also \cite[Prop. 8.3]{GvN}. 

Whenever $\nabla_H$ is closable in $L^p(X, \mu_{\infty})$, the Sobolev space $W^{1,p}_Q(X, \mu_{\infty})$ is defined as the domain of its closure (still called $\nabla_H$), and it is a Banach space endowed with the graph norm
$$
\|f\|_{W^{1,p}_Q(X, \mu_{\infty})}^p = \|f\|_{L^{p}(X, \mu_{\infty})}^p + \int_X \|\nabla_Hf(x)\|_H^p\,\mu_{\infty}(dx) . 
$$ 
In particular, for $p=2$ it is a Hilbert space with inner product $\langle f, g\rangle _ {W^{1,2}_Q(X, \mu_{\infty})} = \langle f, g\rangle _ {L^2(X, \mu_{\infty})}
+ \langle \nabla_H  f,  \nabla_Hg\rangle _ {L^2(X, \mu_{\infty};H)}$. 
In its turn, the operator $D_H^2:  {\mathcal F}_0 \subset L^p(X, \mu_{\infty}) \to L^p(X, \mu_{\infty}; {\mathscr L}_2(H))$   is closable, 
and $W^{2,p}_Q(X, \mu_{\infty})$ is defined as the domain of the closure (still called $D^2_H$),   endowed with the graph norm
$$
\|f\|_{W^{2,p}_Q(X, \mu_{\infty})}^p = \|f\|_{W^{1,p}_Q(X, \mu_{\infty})}^p + 
\int_X \|D^2_Hf(x)\|_{ {\mathscr L}_2(H)} ^p\,\mu_{\infty}(dx) . 
$$
Another involved Sobolev-type space is the domain of the closure of 
$A^*_{\infty}\nabla_{H_{\infty}}:  {\mathcal F}_0 \subset L^p(X, \mu_{\infty}) \to L^p(X, \mu_{\infty}; H_{\infty})$ 
in  $L^p(X, \mu_{\infty})$, called $W^{1,p}_{AQ}(X, \mu_{\infty})$  (we recall that $A_{\infty}$ is the part of $A$ 
in $H_{\infty}$). 

Using the notation in \eqref{F0}, for $f\in {\mathcal F}_0$ we have $\|\nabla_Hf(x)\|_H = \|Q^{1/2 }\nabla f(x)\|$, 
$\|D^2_Hf(x)\|_{{\mathscr L}_2(H)}^2 ={\rm Tr}\, (QD^2f(x))^2$, and 
$\|A^*_{\infty}\nabla_{H_{\infty}}f(x)\|_{H_{\infty}}^2 = 
\langle  A^*\nabla \varphi(x), Q_{\infty}A^*\nabla \varphi(x)\rangle $. 
In the symmetric case, using the Lyapunov equation we get 
$\|A^*_{\infty}\nabla_{H_{\infty}}f(x)\|_{H_{\infty}}^2= \langle \nabla \varphi (x), -AQ   \nabla \varphi(x) \rangle/2$.  
In the case of the classical Ornstein-Uhlenbeck operator, we have $A= -I$, 
$Q_{\infty} = 2Q$, and the spaces $W^{1,p}_Q(X, \mu_{\infty})=W^{1,p}_{AQ}(X, \mu_{\infty})$, $W^{2,p}_Q(X, \mu_{\infty})$ 
considered here coincide respectively with the spaces $W^{1,p}(X, \gamma)$, $W^{2,p} (X, \gamma)$ described in Section 
\ref{SecPrelimin}(iv), with $\gamma =\mathcal N_{0, 2Q}$. 

Before going on, we observe that the quadratic form
$$
\mathcal Q (\varphi, \psi) := \frac{1}{2}
\int_X \langle \nabla_H\varphi (x), \nabla _H\psi (x)\rangle_H \mu_{\infty}(dx), \quad 
\varphi, \;\psi \in W^{1,2}_Q(X, \mu_{\infty})
$$
is closed, and in the symmetric case $-L_2$ is the operator associated with the form $\mathcal Q$ in 
$L^2(X, \mu_{\infty})$, namely
$$
D(L_2) =\{ f\in W^{1,2}_Q(X, \mu_{\infty});\; \exists g\in L^2(X, \mu_{\infty}) \;\text{s.t.} \; 
\mathcal Q (f, \psi) = \langle f, g\rangle_ {L^2(X, \mu_{\infty})}\}, \quad L_2f = -g, 
$$ 
and therefore $D(-L_2)^{1/2} = W^{1,2}_Q(X, \mu_{\infty})$. Even in the nonsymmetric case, recalling that $\mathcal F_0$ is a core for $L_2$, formula \eqref{integrazioneL} yields  $D(L_2)\subset W^{1,2}_Q(X, \mu_{\infty})$ and \eqref{integrazioneL} holds for any $f$, $g\in D(L_2)$. In particular, taking $f=g$ we get 
\begin{equation}
\label{integrazioneLgenerale}
\int_X  Lf(x) \,f(x) \,\mu_{\infty}(dx) = - \frac{1}{2} \int_X \|\nabla_Hf \|_H ^2\,\mu_{\infty}(dx), \quad f\in D(L). 
 \end{equation}
In the analytic case (see condition (ii) of Thm. \ref{Thm:analyticity}) there is a sort of bounded extension of $Q_{\infty}A^*$ to $H$; more precisely, see \cite{MaasVN1}, there exists an operator $B\in \mathscr L (H)$ such that $BQ_{\infty}x = Q_{\infty}A^*x$ for $x\in D(A^*)$, and that satisfies $B+B^* = -I$ in $H$ by the Lyapunov equation. 
Moreover,   $L_pf = \nabla_H^*B\nabla_Hf$, for every $f$ in the core $ \mathcal{F }_0 $. In the symmetric case we have $B=-I/2$, and this statement coincides with  \eqref{divgrad} for the classical Ornstein-Uhlenbeck operator. 

The next theorem follows from  \cite{Ch-MiGol_JFA,Ch-MiGol1,MaasVN2,MaasVN4}, and generalizes an earlier result of \cite{BeppeLincei}.  

\begin{theorem}
\label{th:Chojn_Gold}
In the symmetric  case for every $p\in (1, +\infty)$ we have
 $D(L_p) = W^{2,p}_Q(X, \mu_{\infty})\cap W^{1,p}_{AQ}(X, \mu_{\infty})$,  $D((-L_p)^{1/2})=  W^{1,p}_Q(X, \mu_{\infty})$, 
with equivalence of the respective norms. 
\end{theorem}

The next theorem follows from  \cite{MaasVN3,MaasVN4}. We recall that in the analytic case $e^{tA}$ maps $H$ into itself, and $S_H(t)= e^{tA}_{|H}:H\to H$ is a strongly continuous semigroup. We denote by $A_H$ its infinitesimal generator.  

\begin{theorem}
\label{th:MVN}
In the analytic  case,  the following conditions are equivalent. 
\begin{itemize}
\item[(i)] $D((-L_p)^{1/2})=  W^{1,p}_Q(X, \mu_{\infty})$, with equivalence of the respective norms; 
\item[(ii)] the operator $-A_H$ admits bounded $H^{\infty}$ functional calculus in $H$. 
\end{itemize}
If such equivalent conditions are satisfied, we have
 $D(L_p) = W^{2,p}_Q(X, \mu_{\infty})\cap W^{1,p}_{AQ}(X, \mu_{\infty})$, with equivalence of the respective norms. 
\end{theorem}

Theorem \ref{th:MVN} is a generalization of \ref{th:Chojn_Gold}, since in the symmetric case (i) and (ii) are satisfied. 

In   \cite{MaasVN2} sufficient conditions were given in order that $D(L_p) \subset  W^{2,p}_Q(X, \mu_{\infty}) $ for 
$p\in (1, 2]$, even in the nonanalytic case. 
 
Concerning summability improving, the following hypercontractivity result holds. 

\begin{theorem}
\label{Fuhrman}
Fix $t>0$ and let $1\leq p < q $ be such that 
\begin{equation}
\label{q(t)}
q-1 \leq (p-1)\|Q_{\infty}^{-1/2}e^{tA}Q_{\infty}^{1/2}\|_{\mathscr L(X)}^{-2}. 
\end{equation}
Then $T_p(t)(L^p(X, \mu))\subset L^q(X, \mu)$, and $\|T_p(t)f\|_{ L^q(X, \mu)}\leq \|f\|_{ L^p(X, \mu)}$ for every $f\in L^p(X, \mu)$. 
\end{theorem}
 
The proof is in \cite{Fuhrman} and (in the case that $Q_{\infty}$ is one to one) in \cite{Ch-MiGol2}. Of course, the statement is  meaningful only if \eqref{eq:Fuhrman} is satisfied. As we mentioned before, if \eqref{controllability} holds then \eqref{eq:Fuhrman} holds for every $t>0$. Another simple example is the case that $Q$ commutes with $e^{tA}$  and \eqref{sgrdecay} holds; then $Q_{\infty}^{-1/2}e^{tA}Q_{\infty}^{1/2}= e^{tA}$ and \eqref{eq:Fuhrman} is satisfied for large $t$ if $M>1$, for every $t>0$ if $M=1$, independently of 
the validity of \eqref{controllability}. In particular, if $A = -\omega I$ with $\omega >0$, \eqref{controllability} is not satisfied but  \eqref{eq:Fuhrman} holds for every $t>0$. 

For the classical Ornstein-Uhlenbeck semigroup of Section \ref{sect:OUclassico} condition \eqref{q(t)} coincides with the hypercontractivity property  stated there.  

It is well known, see \cite{Gross,DaGrSi92}, that under appropriate assumptions the hypercontractivity of a semigroup is equivalent to the occurrence of a suitable logarithmic Sobolev inequality. But for general Ornstein-Uhlenbeck semigroups the assumptions of \cite{Gross} are not necessarily satisfied, as shown in \cite{Fuhrman}. In the symmetric case, namely under the conditions of Theorem \ref{thm:CMG1}, they are satisfied, and by \cite[Thm. 4.2]{Ch-MiGol1} for every $\beta >0$  the following conditions are equivalent. 

\begin{itemize}
\item[(i)] $\|Q^{-1/2}e^{tA}Q^{1/2}\|_{\mathscr L(X)} \leq e^{-\beta t}$, for every $t>0$; 
\item[(ii)] $\|Q_{\infty}^{-1/2}e^{tA}Q_{\infty}^{1/2}\|_{\mathscr L(X)} \leq e^{-\beta t}$, for every $t>0$; 
\item[(iii)]  for every $f\in D(L_2)$ we have 
$$\int_X |f(x)|^2 \log (|f(x)|) \mu_{\infty}(dx) \leq \frac{2}{\beta}\langle -L_2 f, f\rangle _{L^2(X, \mu_{\infty})} + \|f\| _{L^2(X, \mu_{\infty})}^2 \log( \|f\|_{L^2(X,\mu_{\infty}) } ), $$
\item[(iv)] $T(t)$ is a contraction from $L^p(X, \mu_{\infty})$ to $L^q(X, \mu_{\infty})$ for every $t>0$, $1\leq p\leq q$ such that 
$q-1\leq (p-1)e^{2\beta t}$. 
\end{itemize}

In \cite{Fuhrman} it was remarked that if \eqref{eq:Fuhrman}  holds for some $t>0$, then there exist $K$, $\nu>0$ such that 
$$
\left\|T_2(t)f -\int_X f(x)\mu_{\infty}(dx)\right\|_{L^2(X, \mu_{\infty})} \leq Ke^{-\nu t}\|f\|_{L^2(X, \mu_{\infty})}, \quad t>0, \; f\in L^2(X, \mu_{\infty}). 
$$
Notice that the operator $\Pi: L^2(X, \mu_{\infty})\to L^2(X, \mu_{\infty})$, $(\Pi f)(x) = \int_X f(x)\mu_{\infty}(dx)$ for a.e. $x\in X$, is just the orthogonal projection on the subspace of constant functions. 

In general, exponential convergence of $T_2(t)f $ to $\Pi f$ is related to the behavior of the semigroup $S_H(t)$. Indeed, if $e^{tA}$ maps $H$ into itself, for every 
$f\in C^1_b(X)$, $t>0$ and $h\in H$ formula \eqref{comm} yields
$$
\frac{\partial T(t)f}{\partial h}(x) = \int_X \langle \nabla f(e^{tA}x +y), e^{tA}h\rangle _X \mu_t(dy)  = \int_X \langle \nabla_Hf(e^{tA}x +y), e^{tA}h\rangle _H \mu_t(dy), 
$$
and therefore, if $\|S_H(t)\|_{\mathscr L(H)}\leq M_1e^{-\beta t}$ for some $M_1$, $\beta>0$, we argue as in Section \ref{sect:OUclassico} and we obtain
\begin{align*}
| \langle \nabla_HT(t)f(x), h\rangle_H |  &=  
\Bigl| \frac{\partial T(t)f}{\partial h}(x) \Bigr|  
\leq   M_1e^{-\beta t}\|h\|_H  \int_X \|\nabla_Hf(e^{tA}x +y)\| _H \,\mu_t(dy)
\\
&\leq  M_1e^{-\beta t}\|h\|_H \Bigl( \int_X \|\nabla_Hf(e^{tA}x +y)\|^2_H \, \mu_t(dy) \Bigr)^{1/2}
\\
&= M_1e^{-\beta t}\|h\|_H(T(t)\Bigl(\|\nabla_Hf\|_H^2)(x)\Bigr)^{1/2}
\end{align*} 
namely, 
\begin{equation}
\label{decad_gradiente}
\|\langle \nabla_HT(t)f(x)\|_H \leq  M_1e^{-\beta t} \Bigl( T(t)(\|\nabla_Hf\|^2)(x)\Bigr)^{1/2}, 
\quad t>0, \; x\in X. 
\end{equation}
Squaring and integrating with respect to $\mu_{\infty}$ we get, for every  $t>0$, 
$$
\int_X\|\nabla_HT(t)\|_H^2  \,d\mu_{\infty}  \leq  M_1^2e^{-2\beta t}  \int_X  T(t)(\|\nabla_Hf\|^2_H)\,d\mu_{\infty}  = M_1^2e^{-2\beta t}  \int_X \|\nabla_Hf\|^2_H \,d\mu_{\infty}. 
$$
In the analytic case this estimate and \eqref{integrazioneLgenerale} allow to obtain a Poincar\'e inequality, 
\begin{equation}
\label{Poincare}
\int_X | f -\Pi f |^2\,d\mu_{\infty}  \leq \frac{M_1^2}{2\beta} \int_X \|\nabla_H  f\|_H^2\,d\mu_{\infty} , \quad f\in W^{1,2}_Q(X, \mu_{\infty})
\end{equation}
by a classical method that seems to go back to \cite{DS} (the proof given in \cite[Prop. 10.5.2]{DPZbrutto} for a particular case works as well in general, using as main ingredients  \eqref{integrazioneLgenerale} and \eqref{decad_gradiente}). 

By the invariance of $\mu_{\infty}$, $T_2(t)$ maps $L^2_0(X, \mu_{\infty}):= (I-\Pi)(L^2(X, \mu_{\infty}))$ into itself.  Moreover, 
\eqref{Poincare}  and \eqref{integrazioneLgenerale} yield
 $\langle L_2f, f\rangle _{L^2(X, \mu_{\infty})} \leq - (\beta/M_1^2)\|f\|_{L^2(X, \mu_{\infty})}^2$ for every $f\in D(L_2)\cap L^2_0(X, \mu_{\infty})$.
By the general theory of semigroups (e.g. \cite[Section IX.8]{Y}), $\| T_2(t)\|_{\mathscr L (L^2_0(X, \mu_{\infty}))} \leq e^{- \beta t/M_1^2}$ for $t>0$,
  and therefore
\begin{equation}
\label{decayL2}
\|T_2(t)f - \Pi f\|_{L^2(X, \mu_{\infty})} \leq  e^{- \beta t/M_1^2}\| f  \|_{L^2(X, \mu_{\infty})}, \quad t>0, \; f\in L^2(X, \mu_{\infty}). 
\end{equation}
If in addition  \eqref{eq:Fuhrman} holds for some $t>0$,   the rate of convergence of $T_p(t)f$ to 
$ \Pi f$ is the same in all spaces $L^p(X, \mu_{\infty})$, $1\leq p<\infty$. 
  Indeed, 
 if  $p>2$ we fix $\tau>0$ such that $T(\tau)$ is a contraction from $ L^2(X, \mu_{\infty})$ to  $L^p(X, \mu_{\infty})$ 
 (such a $\tau$ exists,  since $Q_{\infty}^{-1/2}e^{tA} Q_{\infty}^{1/2}$ is a semigroup, and therefore if \eqref{eq:Fuhrman} holds for some $t>0$   then $\lim_{\tau\to \infty}\|Q_{\infty}^{-1/2}e^{\tau A} Q_{\infty}^{1/2}\|_{\mathscr L(X)} =0$). For every $t\geq \tau$ and  $f\in L^p(X, \mu_{\infty})$  we have
$$
\|T(t)f -  \Pi f\|_{L^p(X, \mu_{\infty})} = \|T(\tau)(T(t-\tau)f -\Pi f)\|_{L^p(X, \mu_{\infty})} 
\leq \| T(t-\tau)f -\Pi f\|_{L^2(X, \mu_{\infty})}
$$ 
 by Theorem \ref{Fuhrman}, and using  \eqref{decayL2} we get
$$
 \|T(t)f -  \Pi f\|_{L^p(X, \mu_{\infty})} \leq  e^{- \beta (t-\tau)/M_1^2}\| f  \|_{L^2(X, \mu_{\infty})} \leq  e^{\beta  \tau/M_1^2}  e^{- \beta t /M_1^2}\| f  \|_{L^p(X, \mu_{\infty})}, \quad t\geq \tau . 
$$
 Similarly, if $p<2$ we fix $\tau >0$ such that $T(\tau)$ is a contraction from $ L^p(X, \mu_{\infty})$ to  $L^2(X, \mu_{\infty})$.  For every $t\geq \tau$  and  $f\in L^p(X, \mu_{\infty})$  we have 
$$
\|T(t)f -  \Pi f\|_{L^p(X, \mu_{\infty})} \leq \|T(t)f -  \Pi f\|_{L^2(X, \mu_{\infty})} =  
\| T(t-\tau)(T(\tau )f -\Pi(T(\tau )f))\|_{L^2(X, \mu_{\infty})}
$$
so that using  \eqref{decayL2} and then Theorem \ref{Fuhrman} we get 
$$
\|T(t)f -  \Pi f\|_{L^p(X, \mu_{\infty})} \leq   e^{- \beta (t-\tau)/M_1^2}\|  T(\tau )f   \|_{L^2(X, \mu_{\infty})} 
\leq  e^{ \beta  \tau/M_1^2} e^{- \beta t /M_1^2}\| f  \|_{L^p(X, \mu_{\infty})},\ t\geq \tau . 
$$

\section{Ornstein-Uhlenbeck semigroups in Banach spaces}

Many of the results of Section 3 have been extended to the case where $X$ is a separable Banach space. In fact, the already mentioned papers \cite{CarDra,GvN,GGVN,MaasVN1,MaasVN2,MaasVN3,MaasVN4,Priola,PvN,Schmuland} deal with the Banach space case. A survey of the state of the art up to 2003 is in \cite{GvN}.  

As in Section \ref{sect:Hilbert},  $Q\in \mathscr L(X^*, X)$ is a symmetric positive operator, and $A:D(A)\subset X\to X$ is the infinitesimal generator of a strongly continuous semigroup $e^{tA}$ in $X$. As in the Hilbert case, the basic assumption of this section is that for every $t>0$ the operator $Q_t$  defined by \eqref{Qt} is the covariance  of a Gaussian measure $\mu_t$, and in this case the Ornstein-Uhlenbeck semigroup $T(t)$ is defined by \eqref{OUHilbert}. 

If $Q$ itself is a covariance and $A= -I$, $T(t)$ is the classical Ornstein-Uhlenbeck semigroup of Sect. \ref{sect:OUclassico}, provided $\gamma$ is the centered Gaussian measure in $X$ with covariance $2Q$. 

As in the Hilbert case, it is the transition semigroup of a stochastic differential equation in $X$, with an appropriate notion of mild solution  (\cite{BvN,vNW}), and it is a contraction semigroup in ${\mathcal B}_b(X)$ that leaves invariant the spaces $C_b(X)$, $BUC(X)$, $C^{\alpha}_b(X)$, $C^{k}_b(X)$, $C^{\alpha + k}_b(X)$ for $\alpha \in (0, 1)$, $k\in \N$. 

The strong-Feller property of $T(t)$ is not easily recognizable as in the Hilbert case. Characterizations and sufficient conditions for $T(t)$ be strong-Feller are in \cite[Sect. 6.1]{GvN}. 

Concerning the behavior of $T(t)$ in $C_b(X)$, it is strongly continuous in the mixed topology,  and the space 
$\mathcal{F }_0$ defined now by 
\begin{align}
\mathcal{F }_0 := \{  f\in C_b(X):\: &f = \varphi(\langle \cdot, a_1\rangle, \ldots, \langle \cdot, a_n\rangle); \; 
\varphi \in C^{2}_b(\R^n), 
\nonumber \\ \label{F0Banach}
&\;n\in \N, \; a_i\in D(A^*),\ A^*D f(\cdot)(\cdot) \in C_b(X)\}
\end{align}
is a core of the generator $L$ of $T(t)$ in the mixed topology, by \cite[Thm.6.6]{GvN}. The domain of $L$ is still given by \eqref{dominio_generatore}, see \cite[Section 6.1]{GvN}. 

The spaces $H:=H_Q$ and $H_t:=H_{Q_t}$ introduced in Sect. \ref{SecPrelimin}(i)   play the role of the spaces $Q^{1/2}(X)$, $Q_t^{1/2}(X)$ in the Hilbert case. We recall that $H_t$ is the Cameron-Martin space of the measure $\mu_t$. 

An important assumption, already mentioned in Sect.  \ref{sect:Hilbert}, is that $e^{tA}$ maps $H$ into itself, and $S_H(t):= e^{tA}_{|H}:H\to H$ is a strongly continuous semigroup in $H$. As in  in Sect.  \ref{sect:Hilbert}, in this case $e^{tA}$ maps $H$ into $H_t$ for every $t>0$, and $\sup_{0<t<1}\|e^{tA}\|_{\mathscr L(H, H_t)} <\infty$, by \cite[Thm. 3.5]{GvN}. As a consequence, $T(t)$ is smoothing along $H$. See  \cite[Sect. 2]{MaasVN4}, and  \cite{LR} for representation formulae and estimates for any order $H$-derivatives of $T(t)f$ when $f\in C_b(X)$. Again, as in the Hilbert case,  Schauder type theorems were proved in  \cite{LR}, that generalize the one stated in Sect. \ref{sect:OUclassico}, and precisely for every $\alpha \in (0, 1)$, $\lambda >0$ and $f\in C^{\alpha}_H(X)$, $R(\lambda, L)f\in C^2_H(X)$ and $D^2_H R(\lambda, L)f\in C^{\alpha}_H(X, \mathcal L^2(H))$. 

Concerning asymptotic behavior and existence of invariant measures,   assumption \eqref{invmeas} is generalized as follows. 
\begin{equation}
\label{invmeasBanach}  
\left\{ \begin{array}{ll}
 (i) &  \forall f\in X^*\; \;\exists \;\text{weak}-\lim_{t\to \infty} Q_tf := Q_{\infty}f,  
\\
(ii)& Q_{\infty}\; \text{is the covariance of a centered Gaussian measure} \; \mu_{\infty}. 
\end{array}\right. 
\end{equation}
Condition (i) is satisfied if \eqref{sgrdecay} holds, in which case a representation formula similar to \eqref{Qinfty} holds, namely 
$Q_{\infty}f=  \int_0^{\infty}  e^{sA}Qe^{sA^*}f\,ds$ for every $f\in X^*$, where now the integral converges as a Pettis integral, see \cite[Sect. 2]{GvN}. 
As in the Hilbert case, if (i) holds the operator $Q_{\infty}$ maps $D(A^*)$ into $D(A)$ and satisfies the Lyapunov equation \eqref{Lyapunov}; moreover (i) holds iff there exists  a symmetric and positive operator $P\in \mathscr L(X^*, X)$ mapping $D(A^*)$ into $D(A)$ such that $PA^*f + APf = -Qf$ for every $f\in D(A^*)$, see \cite[Sect. 4]{GvN}. 

However, establishing whether a given symmetric positive operator is the covariance of a Gaussian measure is not as simple as in the Hilbert case. Necessary and sufficient conditions are in \cite{vNW}. 
If \eqref{invmeasBanach}   holds, denoting by $H_{\infty}: = H_{Q_{\infty}} =$ the Cameron-Martin space of $ \mu_{\infty}$
(as in the Hilbert case), several statements of the previous section are extendable to the Banach setting. In particular: 
\vspace{1mm}

\noindent (a) $e^{tA}$ maps $H_{\infty}$ into itself, and $e^{tA}_{|H_{\infty}}:H_{\infty}\to H_{\infty}$ is a strongly continuous contraction semigroup, still denoted by $S_{\infty}(t)$. Moreover, for any $t>0$ we have $H_{\infty}= H_{t}$   iff $\|S_{\infty}(t)\|_{\mathscr L(H_{\infty})} <1$. 

\vspace{2mm}
\noindent (b) $ \mu_{\infty}$ is an invariant measure of $T(t)$, and the  arguments used in Sections \ref{sect:OUclassico} and \ref{sect:Hilbert} yield that  $T(t)$ extends to a contraction $C_0$-semigroup  $T_p(t)$ in $L^p(X, \mu_{\infty})$, for every $p\in [1, +\infty)$. 

\vspace{2mm}
\noindent (c) Conditions (i) and (iii) of Theorem \ref{thm:CMG1} are still equivalent, see \cite[Thm. 7.4]{GvN}; if they hold $T_p(t)$ is an analytic contraction semigroup in $L^p(X, \mu_{\infty})$ for every $p\in (1, \infty)$. 

\vspace{2mm}
\noindent (d) Conditions (i), (ii), and (iii) of Theorem \ref{Thm:analyticity} are still equivalent, see \cite[Sect. 8]{GvN}; if they hold $T_p(t)$ is an analytic contraction semigroup  in $L^p(X, \mu_{\infty})$ for every $p\in (1, +\infty)$. The optimal angle of analyticity and  the optimal angle for the bounded $H^{\infty}$ calculus of $-L_p$ were determined in \cite{MaasVN1,CarDra}, respectively,  in the present Banach setting. 

\vspace{2mm}
\noindent (e) Theorems \ref{thm:CMG1}  and  \ref{Thm:analyticity} still hold, where the involved Sobolev spaces $W^{1,p}_Q(X, \mu_{\infty})$, $W^{2,p}_Q(X, \mu_{\infty})$, $W^{1,p}_{AQ}(X, \mu_{\infty})$ are defined in a similar way to the Hilbert case.  See \cite{MaasVN2,MaasVN3,MaasVN4}. 

\vspace{5mm}
\noindent {\bf Acknowledgements.} We thank Jan van Neerven for useful discussions. 
The authors are members of GNAMPA-INdAM and they have been partially supported by MIUR through the research project  PRIN 2015233N54. 


\end{document}